\documentclass[11pt]{article}      
\usepackage{graphicx}
\usepackage{latexsym}
\usepackage{amssymb}
\usepackage[mathscr]{eucal}
\usepackage{color}
\usepackage{bbm}
\usepackage{dsfont}
\newcommand{\sub}[1]{{\mbox{\footnotesize $#1$}}}

\newcommand{\Bp}{\begin{proposition}\begin{sl}}
\newcommand{\EP}[1]{\end{sl}\label{proposition:#1}\end{proposition}}

\newcommand{\Div}{\mbox{\rm div}\,}

\newcommand{\curl}{\mbox{\rm curl}\,}

\newcommand{\Int}[2]{{\displaystyle \int_{ #1}^{ #2}}}

\newcommand{\Lim}[1]{{\displaystyle \lim_{ #1}}}

\newcommand{\Sup}[1]{{\displaystyle \sup_{#1}}}

\newcommand{\beea}{\begin{eqnarray}}

\newcommand{\eeea}{\end{eqnarray}}

\newcommand{\bfe}{{\mbox{\boldmath $e$}} }

\newcommand{\0}{{\mbox{\boldmath $0$}} }

\newcommand{\BF}{\begin{footnotesize}}

\newcommand{\EF}{\end{footnotesize}}

\newcommand{\bi}{\begin{itemize}}

\newcommand{\ei}{\end{itemize}}

\newcommand{\ed}{\end{document}}

\newcommand{\be}{\begin{equation}}

\newcommand{\ba}{\begin{array}}

\newcommand{\ea}{\end{array}}

\newcommand{\ee}{\end{equation}}

\newcommand{\eeq}[1]{\label{eq:#1}\end{equation}}

\newcommand{\real}{{\mathbb R}}

\newcommand{\nat}{{\mathbb N}}

\newcommand{\bfomega}{\mbox{\boldmath $\omega$}}

\newcommand{\bfxi}{\mbox{\boldmath $\xi$}}

\newcommand{\bfx}{\mbox{\boldmath $x$}}

\newcommand{\bfy}{\mbox{\boldmath $y$}}

\newcommand{\bfphi}{\mbox{\boldmath $\varphi$}}

\newcommand{\bfPhi}{\mbox{\boldmath $\Phi$}}

\newcommand{\bfv}{{\mbox{\boldmath $v$}} }

\newcommand{\bfu}{{\mbox{\boldmath $u$}} }

\newcommand{\bfw}{{\mbox{\boldmath $w$}} }

\newcommand{\bfH}{{\mbox{\boldmath $H$}} }

\newcommand{\calc}{{\cal C}}

\newcommand{\calh}{{\cal H}}

\newcommand{\cali}{{\cal I}}

\newcommand{\bfV}{{\mbox{\boldmath $V$}} }

\newcommand{\bfg}{{\mbox{\boldmath $g$}} }

\newcommand{\bfn}{{\mbox{\boldmath $n$}} }

\newcommand{\half}{\mbox{$\frac{1}{2}$}}

\def\Bbb R{\real}

\def\bar{\overline}

\newcommand{\ED}{\end{description}}

\def\tag{\renewcommand{\theequation}}

\newcommand{\Br}{\begin{rem}\begin{rm}}

\newcommand{\Er}{\end{rm}\end{remark}}

\newtheorem{lemm}{Lemma}[section]

\newtheorem{theo}{Theorem}[section]

\newtheorem{rem}{Remark}[section]

\newtheorem{coro}{Corollary}[section]

\newtheorem{exe}{\footnotesize{Exercise}}[section]

\newcommand{\Be}{\begin{exe}\begin{footnotesize}\begin{rm}}

\newcommand{\EE}[1]{\end{rm}\end{footnotesize}\label{exe:#1}\end{exe}}

\newcommand{\Bt}{\begin{theo}\begin{sl}}

\newcommand{\Et}{\end{sl}\end{theorem}}

\newcommand{\Bl}{\begin{lemm}\begin{sl}}

\newcommand{\El}{\end{sl}\end{lemma}}

\newtheorem{proposition}{Proposition}[section]

\newcommand{\eqref}[1]{{\rm (\ref{eq:#1})}}

\newcommand{\Bc}{\begin{coro}\begin{sl}}

\newcommand{\Ec}{\end{sl}\end{coro}}

\newcommand{\ET}[1]{\end{sl}\label{theo:#1}\end{theo}}

\newcommand{\EL}[1]{\end{sl}\label{lemm:#1}\end{lemm}}

\newcommand{\theoref}[1]{{\rm Theorem \ref{theo:#1}}}

\newcommand{\ER}[1]{\end{rm}\label{rem:#1}\end{rem}}

\newcommand{\EC}[1]{\end{sl}\label{coro:#1}\end{coro}}

\newcommand{\remref}[1]{{\rm Remark \ref{rem:#1}}}

\newcommand{\lemmref}[1]{{\rm Lemma \ref{lemm:#1}}}

\begin{document}
\title{On the Problem of Steady Bifurcation of a Falling Sphere\\ in a Navier-Stokes Liquid}
\author{Giovanni P Galdi \smallskip \\ { \small Department of Mechanical Engineering and Materials Science}\\ { \small University of Pittsburgh, USA}}
\date{}
\maketitle
\begin{abstract} We study steady bifurcation for the coupled system body-liquid consisting of a sphere  freely falling in a Navier-Stokes liquid under the action of gravity. In particular we show that, under the assumption that for the bifurcating solution the translational velocity of the sphere is parallel to the gravity, bifurcation takes place provided 1 is a simple eigenvalue of a suitable linear operator and the transverslity property holds. Moreover, we also give sufficient conditions for symmetry breaking. 
\end{abstract}
\renewcommand{\theequation}{\arabic{section}.\arabic{equation}}
\setcounter{section}{0}
\section*{Introduction} 
Bifurcation --steady and time-periodic-- is a common phenomenon in fluid mechanics that, over the past several decades, has been the object of numerous  rigorous investigations. Typical and significant examples are the so called Taylor--Couette and B\'enard--Rayleigh problems, for which one can provide a rather complete picture  of  steady-state bifurcation; see, e.g., \cite{ChIo}, \cite{HaIo}, \cite[Sections 72.7--72.9]{Zei4} and the reference therein. 
\par
The fundamental issue that one has to address when investigating this kind of questions is to find the right {\em functional setting} that allows for the use of general abstract results. In this regard,  it should  be pointed out that the steady bifurcation problems previously mentioned, as  many of those collected in classical fluid mechanics literature, regard flow occurring in a {\em bounded} spatial region. In such a case, the problem is  naturally formulated as a nonlinear equation in a suitable Banach space, with the relevant operator being a compact perturbation of the identity. As a consequence, the corresponding linearization possesses a purely discrete spectrum that enables one to provide bifurcation criteria in terms of the spectral properties of the 0 eigenvalue, like  simplicity and transversality \cite[Chapter IV]{Satt}.\par
However, if the liquid occupies a spatial region that is  {\em unbounded} in all directions, as in the flow past an obstacle, the above functional setting is no longer available because of lack of compactness of the nonlinear operator and, in addition, the linearized operator has a large essential spectrum containing 0 \cite{FaNe1}. The general problem of bifurcation from the essential spectrum has been addressed by several authors; see \cite{Stuart} and the references therein. Though of great interest, the theory there developed, however,
does not apply to the flow of a liquid past a body, because it requires a self-adjointness property that is not satisfied in such a case.
\par  
Recently, we have undertaken  a systematic study of bifurcation, in both steady and time-periodic cases, of a Navier--Stokes flow past a  body  \cite{GaRa}, \cite{Gageom}, \cite{Gace}, \cite{GaB}, \cite{GaB1}. In the particular case of steady bifurcation, in order to overcome the issue of the essential spectrum, we found it more appropriate to  formulate the problem in the ``natural" function class where the steady-state Navier--Stokes problem is well set. This class consists of ``graded" (homogeneous) Sobolev spaces, characterized by requiring different summability properties of the various derivatives involved   \cite{Gageom}, \cite{Gace}. In doing so, we have been able to furnish necessary as well as sufficient conditions for the occurrence of steady bifurcation from  a (steady-state) Navier--Stokes flow past a rigid body. 
It should be emphasized that all the results established in \cite{Gageom}, \cite{Gace} are obtained under the assumption that the motion of the rigid body is {\em prescribed}.
\par
Objective of the present paper is to present a steady bifurcation theory  in the case when  the motion of the body is {\em not} prescribed, thus becoming a further unknown.  To our knowledge, this is the first contribution to the rigorous study  of bifurcation in a liquid-solid interaction problem. More specifically, we shall study the problem of a sphere of constant density that is steadily and freely falling in a Navier--Stokes liquid under the action of gravity, $\bfg$. Here, the driving mechanism is represented by the dimensionless buoyancy, $\lambda$ ({\em Galilei number}), that we take to be positive, namely, the sphere is falling, not rising. The basic motion ${\sf s}_0$, of the coupled system consists then of a constant translation of the sphere with velocity $\bfxi$ parallel to $\bfg$, while the liquid executes a corresponding flow of steady-state nature,  when referred to a frame translating with velocity $\bfxi$ as well. Such a flow is rotationally symmetric around the direction of $\bfxi$. There is an abundant literature about this problem from  both experimental and numerical viewpoints (see \cite{Tan,Nak,Jenny,Jenny1,GOB,Dus} and the references therein) showing that there exists a critical value $\lambda_0$ of Galilei number at which the flow of the liquid is no longer rotationally symmetric, and a new flow sets in that only possesses planar symmetry. Our goal is to find a suitable functional setting that allows us to furnish necessary and sufficient conditions for the existence of such $\lambda_0$. Precisely, under the assumption that the bifurcating solution also translates in the direction of $\bfg$, we show that the  problem can be correctly formulated in a suitable Hilbert space, $\calh$,  constituted by functions having a finite Dirichlet integral and reducing to a rigid motion on the boundary; see \eqref{0_0}. The equation for the bifurcating branch is then written as a nonlinear equation in $\calh$, where the relevant operator, ${\sf F}$, suitably defined on a dense set of $\calh$, can be written as the sum of the identity plus a nonlinear operator that, however, unlike the case of a bounded domain, is {\em not} compact; see \eqref{12_1}. Nevertheless, we show that the derivative of ${\sf F}$ at criticality is Fredholm of index 0 (\lemmref{8}) and this allows us to apply general bifurcation results. In this way we prove that a necessary condition for the occurrence of bifurcation at a certain $\lambda_0$, is that 1 is an eigenvalue for a suitable linear operator $\lambda\,{\sf M}=\lambda\,{\sf M}(\lambda)$ evaluated at $\lambda=\lambda_0$; see \eqref{38} and \remref{2.2}. In order to make this assumption meaningful, we show that, in fact, the intersection of the spectrum of ${\sf M}(\lambda)$ with the positive real semi-line is constituted only by eigenvalues with finite multiplicity, for each fixed $\lambda>0$; see \lemmref{10}. We then prove that the above condition is also sufficient, provided 1 is simple and the transversality condition holds; see \theoref{2.1}. It seems quite remarkable that, {\em formally}, we can formulate bifurcation conditions  very similar to those given for flow in a bounded region.   Finally, we address the question of symmetry breaking. As suggested by numerical test \cite{GOB}, this certainly happens if, at the critical value $\lambda=\lambda_0$, the sphere start spinning transversally. Based on this considerations, we furnish a sufficient condition for the transversal component of the angular velocity to be non-zero at criticality, which guarantees  symmetry breaking; see \theoref{2.2}. 
\par
We end this introductory section by pointing out a significant open question. It concerns  the assumption that in the bifurcating solution the sphere translates parallel to the gravity. We need this  hypothesis for merely technical reasons and we believe that it can be removed. However, to date, the proof of the latter seems to be  away from our reach. 
\par
The plan of the paper is as follows. In Section 1, after formulating the problem, we show that the bifurcation problem can be written as a nonlinear operator equation in an appropriate Hilbert space. Successively, in Section 2, we study the relevant properties of the operators involved in the equation. Finally, in Section 3, we give a proof of the main results.

\section{Formulation of the Problem and its Functional Setting} Consider a sphere, $\mathscr S$, of constant density $\rho_{\mathscr S}$ and radius $R$,  freely falling  under the action of gravity in an otherwise quiescent Navier--Stokes liquid, $\mathscr L$. Let $\mathcal F=\{O,\bfe_1,\bfe_2,\bfe_3\}$ be the frame with the origin at the center of $\mathscr S$ and the axis $\bfe_1$ oriented along the acceleration of gravity $\bfg$. The steady-state motions of the coupled system $\mathscr S\cup\mathscr L$ in $\mathcal F$ are then governed by the following set of nondimensional equations \cite[Section 4]{Gah}
\be\ba{cc}\medskip\left.\ba{ll}\medskip\Delta\bfv+\lambda\,\bfxi\cdot\nabla\bfv=\lambda\,\bfv\cdot\nabla\bfv+\nabla p\\
\Div\bfv=0\ea\right\}\ \ \mbox{in $\Omega$}\\ \medskip
\bfv=\bfxi+\bfomega\times\bfx\ \ \mbox{at $\partial\Omega$}\,,\ \ \Lim{\mbox{\footnotesize $|\bfx|$}\to\infty}\bfv(x)=\0\,,\\ 
\Int{\partial\Omega}{}\mathbb T(\bfv,p)\cdot\bfn=\lambda\,\bfe_1\,,\ \   
\Int{\partial\Omega}{}\bfx\times\mathbb T(\bfv,p)\cdot\bfn=0\,.
\ea
\eeq{1}
Here $\bfv,p$  are (non-dimensional) velocity and  pressure fields of $\mathscr L$, while $\bfxi,\bfomega$  represent (non-dimensional) translational and  angular velocities  of $\mathscr S$. Moreover, $\lambda=\sqrt{\alpha\,g\,R^3}/\nu$ ($>0$) is the (dimensionless) Galilei number, where $\alpha=(\rho_{\mathscr S}/\rho_{\mathscr L}-1)$  with $\rho_{\mathscr L}$ density of the liquid,  and $\nu$ its kinematic viscosity. We are assuming that the sphere has a negative buoyancy, so that $\alpha>0$. Finally, $\Omega$ is the region occupied by the liquid (the exterior of $\mathscr S$), and
$$
\mathbb T(\bfv,p)=-p\,\mathds{1}+2\,\mathbb D(\bfv)\,,\ \ \mathbb D(\bfv):=\half\big(\nabla\bfv+(\nabla\bfv)^\top\big)\,,
$$
is the Cauchy tensor with $\mathds{1}$ identity tensor and $\top$ denoting transpose.\par
On physical ground, it is expected  that, for any (positive) value of the Galilei number, \eqref{1}  always admits a solution where $\mathscr S$ moves with no spin and  translational velocity directed along the gravity. In this regard, 
in \cite[Theorem 4.7]{Gah} --completed with the results of \cite[Section X.6]{Gab}--  the following theorem is proved.
\Bt
For any given $\lambda>0$, problem \eqref{1} has at least one solution ${\sf s}_0={\sf s}_0(\lambda):=(\bfv_0,p_0,\bfxi_0,\bfomega_0)(\lambda)$ such that\footnote{We use standard notation: $L^q$ is the Lebesgue space with norm $\|\cdot\|_q$, and  $W^{1,2}$ is Sobolev space. Furthermore, $D^{k,t}$ represent homogeneous Sobolev spaces with semi-norm $\sum_{|l|=k}\|D^lu\|_t$.} 
\be\ba{ll}\medskip
\bfv_0\in D^{2,s}(\Omega)\cap D^{1,r}(\Omega)\cap L^q(\Omega)\cap C^\infty(\Omega),\,;\\ \medskip
p_0\in D^{1,s}(\Omega)\cap L^{s_1}(\Omega)\cap C^\infty(\Omega)\,,\ \ 
\\
\bfxi_0=\xi_0\,\bfe_1\,, \ \xi_0>0\,;\ \ \bfomega_0=\0\,,
\ea
\eeq{000}
{all $s>1, r>4/3, q>2, s_1>3/2$},
and satisfying the ``energy equality"
\be
\|\mathbb D(\bfv_0)\|_2^2=\lambda\,\xi_0\,.
\eeq{WH}
This solution is rotationally symmetric around the $x_1$-axis.
\ET{1}

The main question we want to address is to establish necessary and sufficient conditions for the occurrence of  steady bifurcation for the branch ${\sf s}_0={\sf s}_0(\lambda)$, at some Galilei number $\lambda=\lambda_0$, when the  the translational velocity remains parallel to the gravity direction $\bfe_1$. 
\par
In order to accomplish this goal, we proceed as follows. In the first place, we will equivalently rewrite \eqref{1} as an operator equation in an appropriate function space. Successively, we will show certain fundamental characteristics of the involved operator, notably, its Fredholm property. In this way, we shall finally be able to apply  classical abstract bifurcation theorems and obtain the desired results.
\par
We begin to introduce some function spaces. Let
$$\ba{ll}\medskip
\calc=\mathcal C(\Omega)=\big\{\bfphi\in C_0^\infty(\bar{\Omega}): \Div\bfphi=0 \ \mbox{in $\Omega$}\,;\\ \medskip
\qquad\qquad\qquad\quad\bfphi(x)=\xi_{\mbox{\footnotesize $\bfphi$}}\bfe_1+\bfomega_{\sub{\bfphi}}\times\bfx, \ \xi_{\sub{\bfphi}}\in\real\,,\ \bfomega_{\sub{\bfphi}}\in\real^3\,, \mbox{in a neighborhood of $\partial\Omega$}\big\}\,,\\
\calc_0=\mathcal C_0(\Omega):=\{\bfphi\in\calc:\ \bfphi=\0\,\ \mbox{at $\partial\Omega$}\}\,,
\ea
$$
and  define
$$\ba{ll}\medskip
\mathcal H=\calh(\Omega)\equiv \,\big\{\mbox{completion of $\calc(\Omega)$ in the norm $\|\mathbb D(\cdot)\|_2$}\big\}\,,\\
\mathcal H_0=\calh_0(\Omega)\equiv \,\big\{\mbox{completion of $\calc_0(\Omega)$ in the norm $\|\mathbb D(\cdot)\|_2$}\big\}\,
\ea
$$
The proof of the next lemma is given in \cite[Lemmas 9--11]{Gah}.
\Bl  $\calh$ is a Hilbert space endowed with the scalar product
\be
\big[\bfu,\bfw\big]=\int_\Omega\mathbb D(\bfu):\mathbb D(\bfw)\,,\ \ \bfu,\bfw\in\calh\,,
\eeq{sp}
and the following characterizations hold
\be
\calh=\big\{\bfu\in W^{1,2}_{\rm loc}(\bar{\Omega}: \bfu\in L^6(\Omega),\,\mathbb D(\bfu)\in L^2(\Omega)\,;\ \Div\bfu=0\ \mbox{in\, $\Omega$}\,;\ \bfu(y)=\xi_\sub{\bfu}+\bfomega_\sub{\bfu}\times \bfy\,,\ y\in\partial\Omega\big\}\,.
\eeq{0_0}
and
$$
\mathcal H_0=\big\{\bfu\in\calh:\ \bfu=\0\ \mbox{at $\partial\Omega$}\}.
$$
Moreover, we have
\be
\|\nabla\bfu\|_2\le\sqrt{2}\|\mathbb D(\bfu)\|_2\le 2\|\nabla\bfu\|_2\,,
\eeq{WH1}
and
\be
\|\bfu\|_6\le c_0\,\|\mathbb D(\bfu)\|_2\,,\ \ \bfu\in\calh\,,
\eeq{SoB}
for some numerical constant $c_1>0$.\footnote{Recall that, in our non-dimensionalization, the sphere has radius 1.} Finally, 
there is another positive numerical constant $c_0$ such that
\be
|\xi_\sub{\bfu}|+|\bfomega_\sub{\bfu}|\le c_1\,\|\mathbb D(\bfu)\|_2\,.
\eeq{WH2}
\EL{1}
\par
Let $\calh^{-1}$ be the (strong) dual of $\mathcal H$, and denote by  $\langle\cdot,\cdot\rangle$ and $\|\cdot\|_{-1}$ the corresponding duality pair and associated norm, respectively. We then introduce the following space that will play a fundamental role in our work:\footnote{$\partial_1\equiv\partial/\partial x_1$.}
$$
\mathscr X(\Omega):=\big\{\bfu\in\calh(\Omega): \ \partial_1\bfu\in\calh^{-1}\big\}
$$ 
where $\partial_1\bfu\in\calh^{-1}$ means that there is $C=C(\bfu)>0$ such that
\be
|(\partial_1\bfu,\bfphi)|\le C\,\|\mathbb D(\bfphi)\|_2\,,\ \ \mbox{for all $\bfphi\in\calc$}\,.
\eeq{3}
In fact, if \eqref{3} holds, since $\calc$ is dense in $\calh$, by the Hahn--Banach theorem
$\partial_1\bfu$ can be uniquely extended to a bounded linear functional on the whole of $\calh$, with
$$
\|\partial_1\bfu\|_{-1}:=\sup_{\mbox{$\bfphi\in\calc; \|\mathbb D(\bfphi)\|_2=1$}}|(\partial_1\bfu,\bfphi)|\,.  
$$
It is readily shown that the functional 
$$
\|\bfu\|_{\mathscr X}:=\|\mathbb D(\bfu)\|_2+\|\partial_1\bfu\|_{-1}\,,
$$
defines a norm in $\mathscr X(\Omega)$. We have the following property whose proof is entirely analogous to \cite[Proposition 65]{Gace} and therefore it will be omitted.
\Bl The space
$\mathscr X(\Omega)$ endowed with the norm $\|\cdot\|_{\mathscr X}$ is a reflexive, separable Banach space.
\EL{0}
\smallskip\par
We also have.
\Bl The space $\mathscr X(\Omega)$ is continuously embedded in $L^4(\Omega)$. 
Moreover, there is $c>0$ such that
$$
\|\bfu\|_4\le c\,\big(\|\partial_1\bfu\|_{-1}^{\frac14}\|\mathbb D(\bfu)\|_2^{\frac34}+\|\mathbb D(\bfu)\|_2\big)\,.
$$
\EL{2}
{\em Proof.} For a given $\bfu\in\mathscr X(\Omega)$, we set
$$
\bfV=-\half\curl[\curl(\zeta\,\xi_\sub{\bfu}\bfe_1x_2^2)+\zeta|\bfx|^2\bfomega_\sub{\bfu}]\,,
$$
where $\zeta$ is a function in $C_0^\infty(\bar{\Omega})$ that is 1 near $\partial\Omega$. Clearly, $\bfV$ is smooth with  bounded support and
$$
\Div\bfV=0\,\ \mbox{in $\Omega$}\,,\ \ \bfV(x)=\xi_\sub{\bfu}+\bfomega_\sub{\bfu}\times\bfx\,\ \mbox{at $\partial\Omega$}\,,
$$
so that $\bfV\in\calc$. 
Moreover, by a straightforward calculation and with the help of \lemmref{1}, one shows
\be
\|\bfV\|_4+\|\nabla\bfV\|_2+\|\partial_1\bfV\|_{-1}\le C\,(|\xi_\sub{\bfu}|+|\bfomega_\sub{\bfu}|)\le C_1\|\mathbb D(\bfu)\|_2\,. 
\eeq{4}
We next write 
\be
\bfu=\bfu-\bfV +\bfV:=\bfw+\bfV\,.
\eeq{5} 
and observe that, also by the above properties of $\bfV$, we have 
\be
\bfw\in \calh_0\,.
\eeq{6} 
Furthermore,   
from the identity $\Delta\bfphi=2\,\Div\mathbb D(\bfphi)$, $\bfphi\in \calc_0$, and the density of $\calc_0$ in $\calh_0$ (see \lemmref{1}), we easily obtain
\be
\|\nabla\bfw\|_2=\sqrt{2}\,\|\mathbb D(\bfw)\|_2\,,\ \ \bfw\in\calh_0\,.
\eeq{7}
Thus, in particular, with the help of \eqref{4},   \eqref{5}, we get
\be\ba{rl}\medskip
|\partial_1\bfw|_{-1}&\!\!\!:= \Sup{\mbox{$\bfphi\in\calc_0; \|\nabla\bfphi\|_2=1$}}|(\partial_1\bfw,\bfphi)|\le \mbox{$\frac1{\sqrt{2}}$}\,\Sup{\mbox{$\bfphi\in\calc_0; \|\mathbb D(\bfphi)\|_2=1$}}|(\partial_1\bfw,\bfphi)|\\
&\le \mbox{$\frac1{\sqrt{2}}$}\,\|\partial_1\bfw\|_{-1}\le C\,\left(\|\partial_1\bfu\|_{-1}+\|\mathbb D(\bfu)\|_2\right)\,.
\ea
\eeq{8}
From \eqref{6}--\eqref{8} and \cite[Proposition 1.1]{Gageom}, we then conclude, $\bfw\in L^4(\Omega)$ which, in turn, by \eqref{4} and \eqref{5} implies $\bfu\in L^4(\Omega)$. Furthermore, again by  \cite[Proposition 1.1]{Gageom} we have
$$
\|\bfw\|_4\le c_1\left(|\partial_1\bfw|_{-1}^{\frac14}\|\nabla\bfw\|_2^{\frac34}+\|\nabla\bfw\|_2\right)\,,
$$
that once combined with \eqref{4}, \eqref{5}, \eqref{7} and \eqref{8}  completes the proof of the lemma.
\par\hfill$\square$\par
With the above two lemmas in hand, we are now able to write \eqref{1} as an operator equation in the space $\calh^{-1}$. 
To do this, we introduce the following {\em notation}:
$$
\mbox{For $\bfphi\in\calh$, we denote by $\bar{\bfphi}:=\xi_\sub{\bfu}\bfe_1+\bfomega_\sub{\bfu}\times\bfx$ its trace at $\partial\Omega$.}
$$ 
Consider now the operators
$$\ba{ll}\medskip
{\sf S}:\bfv\in \mathscr X(\Omega)\mapsto {\sf S}\bfv \in \calh^{-1}\,,\\
\medskip
{\sf g}:\bfv\in \mathscr X(\Omega)\mapsto {\sf g} \in \calh^{-1}\,,\\
{\sf N}:\bfv\in \mathscr X(\Omega)\mapsto {\sf N}(\bfv) \in \calh^{-1}\,,
\ea
$$
where
\be\left.\ba{ll}\medskip
\langle {\sf S}\bfv,\bfphi\rangle:=\big(\mathbb D(\bfv),\mathbb D(\bfphi)\big)\,;\\ \medskip
\langle {\sf g},\bfphi\rangle:=\xi_\sub{\bfphi}\,;\\
\langle{\sf N}(\bfv),\bfphi\rangle:=(\bfv\cdot\mathbb D(\bfphi),\bfv)+\xi_\sub{\bfv}\,\langle\partial_1\bfv,\bfphi\rangle\,;
\ea\right.\ \  \bfphi\in\calh\,.
\eeq{12_0}
By using \lemmref{1}, \lemmref{2} and Schwarz inequality, it is at once established that all the operators introduced above are well defined. Furthermore, it is readily checked that problem \eqref{1} can be reformulated as the following operator equation
\be
{\sf S}\bfv-\lambda\,{\sf g}-\lambda\,{\sf N}(\bfv)=0\ \ \mbox{in $\calh^{-1}$}.
\eeq{12}
Actually, if we dot-multiply by $\bfphi\in\calc$ both sides of \eqref{1}$_1$, integrate by parts over $\Omega$ and take into account  \eqref{1}$_{2,3,5,6}$ we get
\be
(\mathbb D(\bfv),\mathbb D(\bfphi))-\lambda\,\xi_\sub{\bfphi}=\lambda(\bfv\cdot\mathbb D(\bfphi),\bfv)-\lambda\int_{\partial\Omega}\bar{\bfv}\cdot\bfn (\bar{\bfv}\cdot\bar{\bfphi})+\lambda\,\xi_\sub{\bfv}\,(\partial_1\bfv,\bfphi)\,,\ \ \bfphi\in\calc\,.
\eeq{13}
Denote by $\cali$ the surface integral in \eqref{13}. Observing that, for $x\in\partial\Omega$, $\bfx$ and $\bfn$ are parallel vectors, and $\int_{\partial\Omega}\bfn=\0$, we get
$$
\cali=\int_{\partial\Omega}\big[\xi_\sub{\bfv}^2\,\bfe_1\cdot\bfn\,\bfe_1\cdot(\bfomega_\sub{\bfphi}\times\bfx)+\xi_\sub{\bfphi}\xi_\sub{\bfv}\,\bfe_1\cdot\bfn\,\bfe_1\cdot(\bfomega_\sub{\bfv}\times\bfx)+\xi_\sub{\bfv}\,\bfe_1\cdot\bfn\,(\bfomega_\sub{\bfphi}\times\bfx)\cdot(\bfomega_\sub{\bfv}\times\bfx)\big]:=\cali_1+\cali_2+\cali_3\,.
$$
Recalling that $\real^3-\bar{\Omega}=\Omega_0$,  by Gauss theorem it follows that
$$
\cali_1=\xi_\sub{\bfv}^2\Int{\Omega_0}{}\bfe_1\cdot\nabla(\bfomega_\sub{\bfphi}\times\bfx)\cdot\bfe_1=\xi_\sub{\bfv}^2\Int{\Omega_0}{}\bfe_1\cdot\mathbb D(\bfomega_\sub{\bfphi}\times\bfx)\cdot\bfe_1=0
$$
and, analogously,
$$
\cali_2=0\,.
$$
Moreover, again by Gauss theorem,
$$
\cali_3=\xi_\sub{\bfv}\int_{\Omega_0}\big[\partial_1(\bfomega_\sub{\bfv}\times\bfx)\cdot(\bfomega_\sub{\bfphi}\times\bfx)+\partial_1(\bfomega_\sub{\bfphi}\times\bfx)\cdot(\bfomega_\sub{\bfv}\times\bfx)\big]=0\,.
$$
Thus, 
\be
\int_{\partial\Omega}\bar{\bfv}\cdot\bfn (\bar{\bfv}\cdot\bar{\bfphi})=0
\eeq{int}
and, as a result, \eqref{13}, by density, gives \eqref{12}. Conversely,  \eqref{12} implies \eqref{13}, and so choosing in \eqref{13} $\bfphi\in\calc_0$ we obtain
$$
(\mathbb D(\bfv),\mathbb D(\bfphi))=\lambda(\bfv\cdot\mathbb D(\bfphi),\bfv)+\lambda\,\xi_\sub{\bfv}\,(\partial_1\bfv,\bfphi)\,,\ \ \bfphi\in\calc\,,
$$
which, after integrating by parts and using the arbitrariness of $\bfphi\in\calc_0$ shows that there exists a suitable pressure field $p$ such that $(\bfv,p,\xi_\sub{\bfv}\,\bfe_1,\bfomega_\sub{\bfv})$ satisfies \eqref{1}$_{1,2,3}$. If we  integrate by parts the first and the third term in \eqref{13}, the latter furnishes
$$
\big(\Div\mathbb T(\bfv,p)+\lambda\partial_1\bfv-\bfv\cdot\nabla\bfv,\bfphi)=\int_{\partial\Omega}\bar{\bfphi}\cdot\mathbb T(\bfv,p)\cdot\bfn-\lambda\,\xi_\sub{\bfphi}\,,\ \mbox{for all $\bfphi\in\calc$.}
$$
However, the left-hand side of this equation vanishes because $\bfv,p$ and $\xi_\sub{\bfv}$ satisfy \eqref{1}$_1$, so that the right-hind side vanishes as well and, by the arbitrariness of $\bfphi$, this implies that also \eqref{1}$_{4,5}$ are satisfied.
\par
Our main goal is to investigate bifurcation of the branch ${\sf s}_0(\lambda):=(\bfv_0,p_0,\bfxi_0,\0)(\lambda)$ obtained in \theoref{1}. In this regard, we premise the following lemma.
\Bl Let ${\sf s}_0(\lambda):=(\bfv_0,p_0,\bfxi_0,\0)(\lambda)$ be the solution in \theoref{1} corresponding to $\lambda>0$. Then $\bfv_0\in\mathscr X(\Omega)$. 
\EL{300}
{\em Proof.} We only need to show
\be
\|\partial_1\bfv_0\|_{-1}<\infty\,.
\eeq{9}
We recall that ${\sf s}_0$ solves the following problem
\be\ba{cc}\medskip\left.\ba{ll}\medskip\Delta\bfv_0+\lambda \,\xi_0\,\partial_1\bfv_0=\lambda\,\bfv_0\cdot\nabla\bfv_0+\nabla p_0\\
\Div\bfv_0=0\ea\right\}\ \ \mbox{in $\Omega$}\\ \medskip
\bfv_0=\xi_0\,\bfe_1\ \ \mbox{at $\partial\Omega$}\,,\ \ \Lim{\mbox{\footnotesize $|\bfx|$}\to\infty}\bfv_0(x)=\0\,,\\ 
\Int{\partial\Omega}{}\mathbb T(\bfv_0,p_0)\cdot\bfn=\lambda\,\bfe_1\,,\ \   
\Int{\partial\Omega}{}\bfx\times\mathbb T(\bfv_0,p_0)\cdot\bfn=0\,.
\ea
\eeq{10}
If we dot-multiply both sides of \eqref{10}$_1$ by $\bfphi\in\calc$, integrate by parts over $\Omega$ and use \eqref{10}$_{2,3,5,6}$ we get
$$
\lambda\,\xi_0\,(\partial_1\bfv_0,\bfphi)=(\mathbb D(\bfv_0),\mathbb D(\bfphi))-\lambda\,\xi_\sub{\bfphi}-\lambda\,(\bfv_0\cdot\mathbb D(\bfphi),\bfv_0)\,.
$$
Thus, using  in the latter Schwarz inequality and \lemmref{1},  we infer
$$
\lambda\,\xi_0\,|(\partial_1\bfv_0,\bfphi)|\le c\,\left(\|\mathbb D(\bfv_0)\|_2+\lambda+\lambda\,\|\bfv_0\|_4^2\right)\|\mathbb D(\bfphi)\|_2\,,
$$
and since, by \theoref{1}, $\bfv_0\in L^4(\Omega)$, the proof of the lemma is completed.
\par\hfill$\square$\par

We thus write in \eqref{12} $\bfv=\bfv_0+\bfu$ to obtain the following equation for $\bfu\in \mathscr X(\Omega)$
\be
{\sf L}_\lambda\bfu-\lambda\,{\sf K}_\sub{\bfv_0}\bfu-\lambda\,{\sf N}(\bfu)=0\ \ \mbox{in $\calh^{-1}$}\,,
\eeq{12_1}
where, for all $\bfphi\in\calh$, 
\be\ba{ll}\medskip
\langle{\sf L}_\lambda\bfu,\bfphi\rangle:=({\sf S}\bfu,\bfphi)-\lambda\,\xi_0\,\langle\partial_1\bfu,\bfphi\rangle\\
\langle{\sf K}_\sub{\bfv_0}\bfu,\bfphi\rangle:=2(\bfu\cdot\mathbb D(\bfphi),\bfv_0)+\xi_\sub{\bfu}\langle\partial_1\bfv_0,\bfphi\rangle\,.
\ea
\eeq{12_2}
\Br Since, in general $\bfv_0$ is a function of the Galilei number $\lambda$, so is ${\sf K}_\sub{\bfv_0}$. Whenever needed, we shall emphasize this property by writing ${\sf K}_\sub{\bfv_0}(\lambda)$.\par\hfill$\triangle$\par
\ER{1}
\par
It is plain  that $(\lambda_0,\bfv_0(\lambda_0))$ is a bifurcation point for \eqref{12_0} if and only if $(\lambda_0,\0)$ is a bifurcation point for \eqref{12_1}. Therefore, our problem reduces to find a  branch of non-trivial solutions $(\lambda,\bfu(\lambda))$ to \eqref{12_1} in a neighborhood of $(\lambda_0,\0)$. To reach this goal, we need several preparatory results concerning the relevant functional properties of the operators defined in \eqref{12_1}. 
\setcounter{equation}{0}
\section{Preparatory Results} 
We begin with the following.
\Bl For any $\bfu\in\mathscr X(\Omega)$ we have
$$
\langle\partial_1\bfu,\bfu\rangle=0\,.
$$
\EL{3_0}
{\em Proof.} From the splitting \eqref{5} and \eqref{4} we get $\bfw, \bfV\in\mathscr X(\Omega)$ and  
\be
\langle\partial_1\bfu,\bfu\rangle=\langle\partial_1\bfw,\bfw\rangle+(\partial_1\bfw,\bfV)+(\partial_1\bfV,\bfw)+(\partial_1\bfV,\bfV)\,.
\eeq{p1}
By integrating by parts and recalling that $\bfw=0$ at $\partial\Omega$ we show
\be
(\partial_1\bfw,\bfV)+(\partial_1\bfV,\bfw)=0\,.
\eeq{p2}
Also, recalling that $\real^3-\bar{\Omega}=\Omega_0$, by another integration by parts we obtain
\be
(\partial_1\bfV,\bfV)=\half\int_{\partial\Omega}\bar{\bfu}\cdot\bar{\bfu}\,{\bar\bfu}\cdot\bfn=\int_{\Omega_0}\bar{\bfu}\cdot\mathbb D(\bar{\bfu})\cdot\bar{\bfu}=0\,.
\eeq{p3}
Thus, from \eqref{p1}--\eqref{p3} it follows that
$$
\langle\partial_1\bfu,\bfu\rangle=\langle\partial_1\bfw,\bfw\rangle\,.
$$
However, $\bfw\in \calh_0$ with $|\partial_1\bfw|_{-1}<\infty$ (see \eqref{8}). Therefore, from \cite[Proposition 1.2]{Gageom} we infer $\langle\partial_1\bfw,\bfw\rangle=0$, which concludes the proof. 

\par\hfill$\square$\par 
\Bl For any fixed $\lambda>0$, the operator
$
{\sf K}_\sub{\bfv_0}={\sf K}_\sub{\bfv_0}(\lambda) $ 
is compact. 
\EL{5}
{\em Proof.} For sake of simplicity, in what follows we will not distinguish between a given sequence and its subsequences. Thus, suppose  $\{\bfu_n\}\subset\mathscr X(\Omega)$ with $\|\bfu_n\|_{\mathscr X}\le M$ and $M$ independent of $n$. By \lemmref{2} this implies, in particular, 
\be
\|\bfu_n\|_4+\|\mathbb D(\bfu_n)\|_2\le M_1\,,
\eeq{15}
for another $M_1$ independent of $n$. Therefore, by \lemmref{2}, there exists $\bfu\in \mathscr X(\Omega)$ such that
$$
\bfu_n \rightharpoonup \bfu \ \mbox{in $\mathscr X(\Omega)$}\,.
$$
Also, from \eqref{WH1}, \eqref{15} and classical compact embedding theorems, we deduce
\be
\bfu_n \rightarrow \bfu \ \mbox{in $L^4(\Omega_R)$, for all $R>R_*$}\,,
\eeq{18}
By \eqref{15} and \eqref{WH2} there exist $\xi_0\in\real$, $\bfomega_0\in\real^3$ such that
\be
 \xi_\sub{\bfu_n}\rightarrow \xi_*\,,\,\ \bfomega_\sub{\bfu_n}\rightarrow \bfomega_*\ \mbox{in $\real$}\,.
\eeq{WH3}
However, by trace theorems (e.g. \cite[Theorem II.4.1]{Gab}) and \eqref{WH1}  we also have
\be
 \int_{\partial\Omega}|\xi_\sub{\bfu_n}-\xi_\sub{\bfu}+ (\bfomega_\sub{\bfu_n}- \bfomega_\sub{\bfu})\times\bfx|^2\le c\,\big(\|\bfu_n-\bfu\|_{4,\Omega_R}^2+\|\bfu_n-\bfu\|_{4,\Omega_R}\|\mathbb D(\bfu_n-\bfu)\|_2\big)\,.
\eeq{17}
Thus, passing to the limit $n\to\infty$ in \eqref{17} and employing \eqref{15}--\eqref{WH3} we infer
\be
\xi_\sub{\bfu_n}\rightarrow \xi_\sub{\bfu}\,,\,\ \bfomega_\sub{\bfu_n}\rightarrow \bfomega_\sub{\bfu}\ \mbox{in $\real$}\,.
\eeq{22}
Consequently,
setting $\bfw_n:=\bfu_n-\bfu$, $\mu_n:=\xi_\sub{\bfu_n}-\xi_\sub{\bfu}$,    from \eqref{12_2} with the help of Schwarz inequality  we prove
$$
\|{\sf K}_\sub{\bfv_0}\bfw_n\|_{-1}\!\le\!\|\bfw_n\|_{4,\Omega_R}\|\bfv_0\|_4
+\|\bfw_n\|_4\|\bfv_0\|_{4,\Omega^R})\\
++|\mu_n|\,\|\partial_1\bfv_0\|_{-1}\,.
$$
If we let $n\to\infty$ into this relation and use \eqref{15}, \eqref{18} and \eqref{22} we get
$$
\limsup_{n\to\infty}\|{\sf K}_\sub{\bfv_0}\bfw_n\|_{-1}\le M_1\,\|\bfv\|_{4,\Omega^R}\,,
$$
which, in turn, by setting $R\to\infty$ and using the absolute continuity of Lebesgue integral, proves the lemma.
\par\hfill$\square$\par
\Bl The operator ${\sf S}\equiv {\sf L}_{0}$ defined in \eqref{12_2}$_1$ is a homeomorphism of $\calh$ onto $\calh^{-1}$, whereas if $\lambda\neq 0$,  ${\sf L}_\lambda$ is a homeomorphism of $\mathscr X(\Omega)$ onto $\calh^{-1}$ .
\EL{7}
{\em Proof.} By Riesz theorem, for any ${\sf f}\in \calh^{-1}$ there is a unique $\bfu\in\calh$ such that 
$$
[\bfu,\bfphi]=\langle {\sf f},\bfphi\rangle\,,\ \mbox{all $\bfphi\in\calh$}\,,
$$
with $[\cdot,\cdot]$ defined in \eqref{sp}, which proves the first part of the lemma. To show the second part, it is enough to show that for any ${\sf f}\in \calh^{-1}$ there is a unique $\bfu\in\calh$ verifying
\be
(\mathbb D(\bfu),\mathbb D(\bfphi))-\lambda\,\xi_0\,(\partial_1\bfu,\bfphi)=\langle{\sf f},\bfphi\rangle\,,\ \ \mbox{for all $\bfphi\in\calc$}
\eeq{29}
and such that
\be
\|\mathbb D(\bfu)\|_2\le \|{\sf f}\|_{-1}\,.
\eeq{30}
In fact, from \eqref{29}, \eqref{30} it also follows
$$
\lambda\,\xi_0\, |(\partial_1\bfu,\bfphi)|\le 2\|{\sf f}\|_{-1}\|\mathbb D(\bfphi)\|_2, \ \ \mbox{for all $\bfphi\in\calc$}\,,
$$
which proves  $\bfu\in\mathscr X(\Omega)$.
Let $\{\bfphi_k\}\subset \calc$ be a basis in $\calh$ with $[\bfphi_k,\bfphi_{k'}]=\delta_{kk'}$.
In view of \eqref{SoB}, we have 
that for any given $\bfphi\in\calh$, there is $\{\gamma_{k}\}\subset \real$ such that setting 
$\bfPhi_N=\sum_{k=1}^N\gamma_k\bfphi_k$, it follows
\be
\|\mathbb D(\bfPhi_N -\bfphi)\|_2+\|\bfPhi_N -\bfphi\|_6\to 0 \ \mbox{as 
 $N\to\infty$}\,.
\eeq{31}
We look for an ``approximating" solution to \eqref{29} of the type
$$
\bfu_m=\sum_{\ell=1}^mc_{\ell m}\bfphi_\ell
$$
where the coefficients $c_{\ell m}$ are solutions to the system
\be
(\mathbb D(\bfu_m),\mathbb D(\bfphi_k))-\lambda\,\xi_0\,(\partial_1\bfu_m,\bfphi_k)=\langle{\sf f},\bfphi_k\rangle\,,\ \ k=1,\ldots,m\,,
\eeq{32}
or, equivalently,
\be
\sum_{k=1}^n\left(\delta_{\ell k}-\lambda \,\xi_0\, \mathcal A_{\ell k}\right)c_{\ell m}=\langle{\sf f},\bfphi_k\rangle\,,\ \ k=1,\ldots,m\,,
\eeq{33}
where
$$
\mathcal A_{\ell k}=(\partial_1\bfphi_{\ell},\bfphi_k)\,.
$$
Integrating by parts, we have
$$
\mathcal A_{\ell k}=-\mathcal A_{k \ell}+\int_{\partial\Omega}\bfe_1\cdot\bfn \,\bar{\bfphi_\ell}\cdot\bar{\bfphi_k}
$$
However, arguing as in the proof of \eqref{int}, we show 
\be
\int_{\partial\Omega}\bfe_1\cdot\bfn \,\bar{\bfphi_\ell}\cdot\bar{\bfphi_k}=0,
\eeq{35}
which implies that $\mathcal A_{\ell k}$ is skew-symmetric. As a result \cite[Lemma IX.3.1]{Gab}, \eqref{33} has a unique solution $c_{\ell k}$, $k=1,\ldots,m$, or, equivalently, \eqref{32} has a unique solution $\bfu_m$ for each $m\in\nat$. Multiplying both sides of \eqref{32} by $c_{\ell k}$,  summing over $k$ from 1 to $m$ and using \eqref{35} entails
$$
\|\mathbb D(\bfu_m)\|_2^2=\langle{\sf f},\bfu_m\rangle
$$
from which we easily deduce
\be
\|\mathbb D(\bfu_m)\|_2\le\|{\sf f}\|_{-1}\,.
\eeq{36}
Thus, there is $\bfu\in\calh$ such that $\bfu_m\rightharpoonup \bfu$ in $\calh$ and, moreover, $\bfu$ satisfies \eqref{30}. Also, passing to the limit $m\to\infty$ in \eqref{32} we show
\be
(\mathbb D(\bfu),\mathbb D(\bfphi_k))-\lambda\,\xi_0\,(\partial_1\bfu,\bfphi_k)=\langle{\sf f},\bfphi_k\rangle\,,\ \ \mbox{for all $k\in\nat$}\,.
\eeq{37}
We can now replace in \eqref{37} $\bfphi_k$ with linear combinations of the type $\bfPhi_N$ given in \eqref{31}. Therefore, taking the limit $N\to\infty$ we obtain that actually $\bfu$ satisfies \eqref{29}, which completes the proof of existence. Suppose, next, that
$$
\langle{\sf L}_\lambda\bfu,\bfphi\rangle\equiv (\mathbb D(\bfu),\mathbb D(\bfphi))-\lambda\,\xi_0\,\langle\partial_1\bfu,\bfphi\rangle=0\,,\ \mbox{for all $\bfphi\in\calh$}\,. 
$$
If we set $\bfphi=\bfu$ into this relation and use \lemmref{3_0} we get $\|\mathbb D(\bfu)\|_2=0$, namely, $\bfu=0$, and uniqueness follows. 
\par\hfill$\square$\par 
From \lemmref{5} and \lemmref{7}  we deduce, in particular, the following property.
\Bl For any fixed $\rho\neq0$, and $\lambda>0$ 
the operator ${\sf L}_\rho-\rho\,{\sf K}_\sub{\bfv_0}(\lambda)$ is Fredholm of index 0.
\EL{8}
\par
We next recall that, by \lemmref{7}, ${\sf S}$ is an homeomorphism from $\calh$ onto $\calh^{-1}$. We may then introduce the operator
\be
{\sf M}:\bfu\in D({\sf M})\equiv\mathscr X(\Omega)\subset \calh \mapsto \ {\sf S}^{-1}[\xi_0\,\partial_1\bfu-{\sf K}_\sub{\bfv_0}\bfu]\in\calh \,.  
\eeq{38}
We recall that, in general, the operator ${\sf M}$ is a function of the Galilei number $\lambda$; see \remref{1}. In the next lemmas we shall show some relevant properties of ${\sf M}$.
\Bl
For each fixed $\lambda>0$, the operator ${\sf M}={\sf M}(\lambda)$ is densely defined and closed.  
\EL{9}
{\em Proof.} Since $\calc\subset \mathscr X(\Omega)\subset \calh$, the density property is obvious. Let $\{\bfv_k\}\subset \mathscr X(\Omega)$, ${\sf m}_k:={\sf M}\bfv_k$, with $\bfv_k\to\bfv$ and ${\sf m}_k\to\bfu$ in $\calh$, for some $\bfv,\bfu\in\calh$. Let us prove that, in fact, $\bfv_k\to\bfv$  in $\mathscr X(\Omega)$. Recalling the definition of ${\sf S}$, we obtain that the relation ${\sf M}(\bfv_k-\bfv_{k'})={\sf m}_k-{\sf m}_{k'}$ is equivalent to
\be
\xi_0\langle\partial_1(\bfv_k-\bfv_{k'}),\bfphi\rangle-2\big(\bfv_0\cdot\mathbb D(\bfphi),(\bfv_k-\bfv_{k'})\big)-(\xi_k-\xi_{k'})\langle\partial_1\bfv_0,\bfphi\rangle=\big(\mathbb D({\sf m}_k-{\sf m}_{k'}),\mathbb D(\bfphi)\big)\,.
\eeq{39}
for all $\bfphi\in\calh$. 
Employing H\"older inequality on the right hand side of \eqref{39},  we show 
\be
\xi_0\,\|\partial_1(\bfv_k-\bfv_{k'})\|_{-1}\le 2\|\bfv_0\|_3\|\bfv_k-\bfv_{k'}\|_6+|\xi_k-\xi_{k'}|\|\partial\bfv_0\|_{-1}+\|\mathbb D({\sf m}_k-{\sf m}_{k'})\|_2\,.
\eeq{40}
Since $(\bfv_k-\bfv_{k'})\to 0$ in $\calh$, by \lemmref{1} we know that
$$
(\bfv_k-\bfv_{k'})\to 0 \ \mbox{in $L^6(\Omega)$}\,;\ \ (\xi_k-\xi_{k'})\to 0 \ \mbox{in $\real$}\,,
$$
and also, from \theoref{1} and \lemmref{300}, $\|\bfv_0\|_3+\|\partial\bfv_0\|_{-1}<\infty$. Therefore, letting $k,k'\to\infty$ it follows  that $\bfv\in\mathscr X(\Omega)\equiv D({\sf M})$, and  $\bfv_k\to\bfv$ in $\mathscr X(\Omega)$. In view of the latter, by formally setting in \eqref{4} $\bfv_{k'}\equiv\bfv$, ${\sf m}_{k'}\equiv\bfu$ and then passing to the limit $k\to\infty$, we conclude ${\sf M}\bfv=\bfu$ which completes the proof that ${\sf M}$ is closed. 
\par\hfill$\square$\par
\Bl For any fixed $\lambda>0$ and $\mu\neq0$ the operator
$$
{\sf H}_\mu=\mu\,{\sf I}-\lambda\,{\sf M}(\lambda)
$$
is Fredholm of index 0.
Furthermore, denoting by $\calh_c$ the complexification of $\calh$, by ${\sf M}_c$ the natural extension of ${\sf M}$ to $\calh_c$ and by $\sigma({\sf M}_c)$ the spectrum of ${\sf M}_c$, we have that $\sigma({\sf M}_c)\cap (0,\infty)$ consists at most of a countable number of eigenvalues of finite algebraic multiplicity that can only cluster at 0.
\EL{10}
{\em Proof.}
Since \be{\sf H}_\mu=\mu\,{\sf S}^{-1}({\sf L}_{\frac1\mu}-\frac1\mu\,{\sf K}_{\bfv_0}):={\sf S}^{-1}\,{\sf T}_\mu\,,\eeq{41} and, by \lemmref{8}, ${\sf T}_\mu$ is Fredholm of index 0, we have
$$
{\rm dim}\,N[{\sf H}_\mu]={\rm dim}\, N[{\sf T}_\mu]=m<\infty\,.
$$
Moreover, from
$$
\calh^{-1}=R({\sf T}_\mu)\oplus S_m
$$
with $S_m$  $m$-dimensional subspace, we deduce that for every $\bfw\in \calh$, it is ${\sf S}\bfw=\bfw_1+\bfw_2$, $\bfw_1\in R({\sf T}_\mu)$, $\bfw_2\in S_m$. Therefore, $\bfw={\sf S}^{-1}\bfw_1 + {\sf S}^{-1}\bfw_2$, with ${\sf S}^{-1}\bfw_1\in R({\sf H}_\mu)$, and ${\sf S}^{-1}
\bfw_2\in {\sf S}^{-1}S_m$, which completes the proof. Clearly, the natural complexification ${\sf H}_{\mu c}$, of ${\sf H}_\mu$ is also Fredholm of index 0, for all $\mu>0$. It then follows that the essential spectrum $\sigma_{\rm ess}({\sf M}_c)$, of ${\sf M}_{c}$ defined as the set of $\mu$ where ${\sf H}_{\mu c}$ is not Fredholm has empty intersection with $(0,\infty)$. We shall next show that the resolvent set $P({\sf M}_c)$ of ${\sf M}_c$ has a non-empty intersection with $(0,\infty)$. Since ${\sf H}_{\mu}$ is Fredholm of index 0, it is enough to show that, for sufficiently large $\mu>0$, it is $N[{\sf H}_\mu]=\{0\}$. From \eqref{41}, we see that the latter is equivalent to 
$$
\mu\,\big(\mathbb D(\bfu),\mathbb D(\bfphi)\big)-\xi_0\,\langle \partial_1\bfu,\bfphi\rangle-2(\bfu\cdot\mathbb D(\bfphi),\bfv_0)-\xi_\sub{\bfu}\langle\partial_1\bfv_0,\bfphi\rangle=0\,,\ \ \bfphi\in\calh.
$$
Choosing $\bfphi=\bfu$ in this relation and using \lemmref{3_0}, \eqref{SoB} and \eqref{WH2} along with H\"older inequality furnishes
$$
\mu\,\|\mathbb D(\bfu)\|_2^2\le C\,\big(\|\bfv_0\|_3+\|\partial_1\bfv_0\|_{-1}\big)\|\mathbb D(\bfu)\|_2^2\,. 
$$
Thus, if $\mu>C\,\big(\|\bfv_0\|_3+\|\partial_1\bfv_0\|_{-1}\big)\equiv\bar{\mu}$, we conclude $\bfu\equiv\0$, namely, $P({\sf M}_c)\cap (\bar{\mu},\infty)\neq\emptyset$. Summarizing, we have shown that $\sigma_{\rm ess}({\sf M}_c)\cap (0,\infty)=\emptyset$ while $P({\sf M}_c)\cap (\bar{\mu},\infty)\neq\emptyset$. Therefore, the stated property about eigenvalues is a consequence of classical results in spectral theory \cite[Theorem XVII.2]{GoGoKa}.\par\hfill$\square$
\par
We conclude this section by observing that, in view of the homemomrphism property of ${\sf S}$, the equation \eqref{12_1} is equivalent to the following one 
\be
{\sf F}(\lambda,\bfu):=\bfu-\lambda\,{\sf M}(\lambda)\bfu-\lambda\,{\sf B}(\bfu)=0\ \ \mbox{in $\calh$\,.}
\eeq{2.1}
with ${\sf B}:={\sf S}^{-1}{\sf N}$.
\setcounter{equation}{0}
\section{Necessary and Sufficient Conditions for Steady Bifurcation.}
The main objective of this section is to investigate steady bifurcation of the flow branch determined in \theoref{1} around Galilei number $\lambda_0$. This  leads us to study, in a neighborhood of $\lambda_0$,  the existence of a nontrivial branch of solutions $\bfu=\bfu(\lambda)$ to equation \eqref{2.1}. 
 Since ${\sf B}$ is quadratic in $\bfu$, we deduce that ${\sf F}$ is analytic with respect to the $\bfu$-variable and, in particular, its Frechet derivative at $(\lambda=\lambda_0,\bfu=\0)$ is given by 
\be
D_\sub{\bfu}{\sf F}(\lambda_0,\0)\bfw=\bfw-\lambda_0\,{\sf M}(\lambda_0)\bfw\,.
\eeq{2_0}
As for the regularity in $\lambda$, we  
 make the following assumption,
\tag{H}
\be 
\mbox{There is a neighborhood $U_0$ of $\lambda_0$ such that the map\ $\lambda\in U_0\mapsto \bfv_0(\lambda)$ is of class $C^2$}
\eeq{H}
which implies that ${\sf F}$ is of class $C^2$ in $U_0\times \mathscr X(\Omega)$. We next observe that, by \lemmref{9}, $\mu\,{\sf I} -\lambda\,{\sf M}(\lambda)$ is Fredholm of index 0 for all fixed $\lambda>0$ and all $\mu\neq 0$. Thus, we may define  a simple eigenvalue of ${\sf M}$ as follows \cite[Definition 79.14]{Zei4}. The number $\mu\neq 0$ is a {\em simple eigenvalue} if
\setcounter{equation}{0}
\renewcommand{\theequation}{\arabic{section}.\arabic{equation}}
\be\ba{ll}\medskip
{\dim}\,N[\mu{\sf I}-\lambda\,{\sf M}(\lambda)]=1\,;
\\
N[\mu{\sf I}-\lambda\,{\sf M}(\lambda)]\cap R[{\sf I}-\lambda\,{\sf M}(\lambda)]=\{0\}\,.
\ea
\eeq{2.2}
It is very well known that the second condition can be reformulated in an equivalent way in terms of an eigenvector of the adjoint operator, ${\sf M}^*$, of ${\sf M}$. Actually, from \eqref{2.2}$_1$ and the Fredholm property we deduce that ${\rm dim}\,N[\mu\,{\sf I}-\lambda\,{\sf M}^*(\lambda)]={\rm codim}\,R[\mu\,{\sf I}-\lambda\,{\sf M}(\lambda)]=1$. Thus, denoting by $\bfw_1\in \calh$ and $\bfw_1^*\in\calh^{-1}$ non-zero elements of $N[\mu\,{\sf I}-\lambda\,{\sf M}(\lambda)]$ and $N[\mu\,{\sf I}-\lambda\,{\sf M}^*(\lambda)]$, respectively, \eqref{2.2}$_2$ is shown to be equivalent (after suitable normalization) to
\be
\langle \bfw_1^*,\bfw_1\rangle=1\,.
\eeq{2.3}
The following result is a consequence  of \eqref{2.1} and \cite[Corollary 79.16]{Zei4}.
\Bl Suppose  there exists $\lambda_0>0$ such that 1 is a simple eigenvalue of the operator $\lambda_0\,{\sf M}(\lambda_0)$ and that \eqref{H} holds. Then, there is  $\bar{U}_0\subseteq U_0$ such that the eigenvalue $\mu=\mu(\lambda)$ of $\lambda\,{\sf M}(\lambda)$, $\lambda\in\bar{U}_0$, is still simple and of class $C^2$. Moreover, we have
$$
\mu'(\lambda_0)=-\langle \bfw_1^*, ({\sf M}(\lambda_0)+\lambda_0\,{\sf M}^\prime(\lambda_0)\big)\bfw_1\rangle\,,
$$
where the prime denotes  differentiation with respect to $\lambda$.
\EL{2.1}
\par
We are now in a position to prove our main bifurcation result.
\Bt A necessary condition for $(\lambda_0,\0)$ to be a bifurcation point of \eqref{2.1} is that $\dim N[{\sf I}-\lambda_0\,{\sf M}(\lambda_0)]\ge1$, namely, the equation
\be
\bfw-\lambda_0\,{\sf M}(\lambda_0)\bfw=\0
\eeq{2.4}
has at least one non-trivial solution $\bfw_1$. Conversely, suppose that $1$ is a simple eigenvalue of $\lambda_0\,{\sf M}(\lambda_0)$, namely, \eqref{2.1} holds with $\mu=1$. Then, if $\mu'(\lambda_0)\neq 0$  (transversality condition), in a suitable neighborhood of $(\lambda_0,\0)$ there exists exactly one continuous curve of nontrivial solutions to \eqref{2.1}, $(\lambda,\bfu(\lambda))$, with $(\lambda_0,\bfu(\lambda_0))=(\lambda_0,\0)$. 
\ET{2.1}{\em Proof.} The necessary condition for bifurcation at $(\lambda_0,\0)$ is that the derivative $D_\sub{\bfu}{\sf F}(\lambda_0,\0)$ is singular. By \eqref{2_0} and \lemmref{10} this derivative is Fredholm of index 0, and, therefore, it is singular if and only if \eqref{2.4} has a nontrivial solution $\bfw_1$. Conversely, since $D_\sub{\bfu}{\sf F}(\lambda_0,\0)$ is Fredholm of index 0, if $\dim N[{\sf I}-\lambda_0\,{\sf M}(\lambda_0)]=1$, a classical bifurcation result \cite[Theorem 4.1.12]{Berger} ensures the stated sufficient property provided
$$ 
D^2_{\lambda\,\mbox{\footnotesize$\bfu$}}{\sf F}(\lambda_0,\0)\,\bfw_1\not\in R[D_\sub{\bfu}{\sf F}(\lambda_0,\0)]\,,
$$
or, equivalently, 
\be
\langle \bfw_1^*,D^2_{\lambda\,\mbox{\footnotesize$\bfu$}}{\sf F}(\lambda_0,\0)\bfw_1\rangle\neq 0\,.
\eeq{2.5}
By a straightforward computation, from \eqref{2.1} we show that
$$
D^2_{\lambda\,\mbox{\footnotesize$\bfu$}}{\sf F}(\lambda_0,\0)\bfw_1={\sf M}(\lambda_0)\bfw_1+\lambda_0\,{\sf M}^\prime(\lambda_0)\bfw_1\,,
$$
so that, if 1 is a simple eigenvalue, by \lemmref{2.1}, condition \eqref{2.5} is equivalent to $\mu'(\lambda_0)\neq 0$, which concludes the proof of the theorem.
\par\hfill$\square$\par
\par
We would like to make several  comments regarding \theoref{2.1} that we shall collect in as many remarks.
\Br
If the branch $\bfv_0(\lambda)$ is constant around $\lambda_0$, the hypothesis $\mu'(\lambda_0)\neq 0$ is equivalent to the request that 1 is a simple eigenvalue and, therefore, the latter can be omitted. 
\par\hfill$\triangle$\par
\ER{2.1}
\Br Taking into account the definition of the operator ${\sf M}$ given in \eqref{38}, we show that \eqref{2.4} is equivalent to the following one
\be
\frac1{\lambda_0}\big(\mathbb D(\bfw),\mathbb D(\bfphi)\big)-\xi_0\,\langle \partial_1\bfw,\bfphi\rangle-2(\bfw\cdot\mathbb D(\bfphi),\bfv_0)-\xi_\sub{\bfw}\langle\partial_1\bfv_0(\lambda_0),\bfphi\rangle=0\,,\ \ \mbox{for all $\bfphi\in\calh$\,,}
\eeq{2.6}
where $\bfw\in\mathscr X(\Omega)$.
Thus, by classical regularity results \cite[Section VII.1]{Gab} and \theoref{1} it follows that, on the one hand, $\bfw\in C^\infty(\Omega)$ and, on the other hand, there is a pressure field $p\in C^\infty(\Omega)$ such that equation \eqref{2.6} is equivalent to the following ones
\be\ba{cc}\medskip\left.\ba{ll}\medskip\Delta\bfw+\lambda_0\,\xi_0\,\partial_1\bfw-\lambda_0\,\big[\bfv_0(\lambda_0)\cdot\nabla\bfw+\bfw\cdot\nabla\bfv_0(\lambda_0)-\xi_\sub{\bfw}\,\partial_1\bfv_0(\lambda_0)\big]=\nabla p\\
\Div\bfw=0\ea\right\}\ \ \mbox{in $\Omega$}\\ \medskip
\bfw=\xi_\sub{\bfw}\,\bfe_1+\bfomega_\sub{\bfw}\times\bfx\ \ \mbox{at $\partial\Omega$}\,,\\ 
\Int{\partial\Omega}{}\mathbb T(\bfw,p)\cdot\bfn=\0\,,\ \   
\Int{\partial\Omega}{}\bfx\times\mathbb T(\bfw,p)\cdot\bfn=\0\,.
\ea
\eeq{2.7}
We omit the proof of the latter, since it is obtained from \eqref{2.6} by an argument entirely analogous to show that \eqref{13} is equivalent to \eqref{1}. In view of the summability properties  of $\bfv_0$ stated in \theoref{1} and the fact that $\bfw\in L^4(\Omega)$, it can be shown that the quantity in bracket  in \eqref{2.6}$_1$ is in $L^t(\Omega)$, for all $t>1$. As a consequence, from known results on the Oseen problem \cite[Section VII.7]{Gab} it follows that $(\bfw,p)$ belongs to the same functional class as $(\bfv_0,p_0)$ given in \eqref{000}. \par\hfill$\triangle$\par  
\ER{2.2}
\Br The bifurcation result given in \theoref{2.1} is of {\em local} nature, namely, the existence of the bifurcating branch is ensured only in a neighborhood of $(\lambda_0,\0)$. As a matter of fact, thanks to \cite[Theorem 7.2]{FPR}, the same  assumptions as those of \theoref{2.1} lead to the following  result of {\em global} nature (see also \cite[Corollary 6.1]{GaRa}).
\smallskip\par\noindent
{\bf Theorem 2.1$^\prime$}\,\ {\sl Let the assumptions of \theoref{2.1} be satisfied and take in \eqref{H}  ${U}_0\equiv(a,b)$. Moreover, denote by $\mathscr C$ the connected component of the closure (in $\real^3\times\mathscr X(\Omega)$) of ${\sf F}^{-1}(0)\backslash (a,b)\times \{\0\}$ that contains $(\lambda_0,\0)$. Then, one of the following three conditions hold:
\begin{itemize}
\item[\rm (i)] $\mathscr C$ contains a point $(a,\bfu)$ or $(b,\bfu)$ for some $\bfu\in\mathscr X(\Omega)$\,;
\item[\rm (ii)] $\mathscr C$ is not compact\,;
\item[\rm (iii)] $\mathscr C$ contains a point $(\lambda_*,\0)$ with $\lambda_*\neq\lambda_0$\,.
\end{itemize}}
It should be observed that, for the validity of the above result, the $C^2$ assumption in \eqref{H} can be replaced by requiring only $C^1$ regularity \cite[Theorem 6.1]{PR}. We also notice that the statement in (ii) could be replaced by ``{\sl $\mathscr C$ is unbounded}," on condition that ${\sf F}$ is proper on the closed bounded sets of $[a+\varepsilon, b-\varepsilon]\times \mathscr X(\Omega)$, for small $\varepsilon>0$. However, the truthfulness of this property does not seem obvious and it is  yet to be ascertained.\par\hfill$\triangle$\par
\ER{2.3}
\par
As already noticed in the introductory section, lab and numerical tests  show that   steady bifurcation occurs by breaking the rotational symmetry of the flow while, however, still keeping planar symmetry along the direction of fall. In particular, the detailed numerical investigation carried out in \cite[Section 3]{GOB} evidences that  symmetry breaking is induced by transversal rotation of the sphere, namely, along a direction perpendicular to the translational velocity. Conversely, it is clear that a non-zero value of such rotation is incompatible with rotational symmetry. Our next objective is then to provide sufficient conditions for the non-vanishing of the transversal component of the angular velocity. Of course, it is enough to furnish such conditions on the solution to \eqref{2.2} or, equivalently, \eqref{2.7}. To this end, let $(\bfH,P)$ solve the following Stokes problem:
\be\ba{cc}\medskip\left.\ba{ll}\medskip
\Div\mathbb T(\bfH,P)=\0\\
\Div\bfH=0\ea\right\}\ \ \mbox{in $\Omega$}\\
\bfH(y)=\bfe_3\times\bfy\,,\ \ \ y\in\partial\Omega\,.
\ea
\eeq{A1}
The fields $\bfH$ and $P$ are well known, and  given by \cite[\S\,334]{Lamb}
\be
\bfH=\bfe_3\times\frac{\bfx}{|x|^3}\,,\ \ P=\textrm{const.}
\eeq{A0}
Clearly, $\bfH\in\calh$, and so we can choose $\bfphi=\bfH$ in \eqref{2.6} to obtain
\be
\frac1{\lambda_0}\big(\mathbb D(\bfw),\mathbb D(\bfH)\big)-\xi_0\,( \partial_1\bfw,\bfH)-2(\bfw\cdot\mathbb D(\bfH),\bfv_0)-\xi_\sub{\bfw}(\partial_1\bfv_0(\lambda_0),\bfH)=0\,.
\eeq{A2}
Notice that we have changed $\langle\cdot,\cdot\rangle$ into $(\cdot,\cdot)$ for the terms involving $\partial_1$, because now they become meaningful in view of the summability properties of $\bfw$ (see \remref{2.2}), $\bfv_0$ (see \theoref{1}) and $\bfH$ (see \eqref{A0}).  
Also, by dot-multiplying both sides of \eqref{A1}$_1$ by $\bfw$,  integrating by parts over $\Omega$ and taking into account \eqref{2.7}$_3$, we deduce
\be
\bfomega_\sub{\bfw}\cdot\int_{\partial\Omega}\bfx\times \mathbb T(\bfH,P)\cdot\bfn=(\mathbb D(\bfw),\mathbb D(\bfH))\,,
\eeq{A3}
where we have used the fact that, as  is well known \cite[p. 187]{HB} and easily  checked,
$$
\int_{\partial\Omega}\bfe_1\cdot\mathbb T(\bfH,P)\cdot\bfn=0\,.
$$
Furthermore,
\be
\int_{\partial\Omega}\bfx\times \mathbb T(\bfH,P)\cdot\bfn=-8\pi\,\bfe_3\,.
\eeq{A4}
Collecting \eqref{A2}--\eqref{A4}, 
we conclude
$$
-\frac8\pi{\lambda_0}\,\bfomega_\sub{\bfw}\cdot\bfe_3=\xi_0\,( \partial_1\bfw,\bfH)+2(\bfw\cdot\mathbb D(\bfH),\bfv_0)+\xi_\sub{\bfw}(\partial_1\bfv_0(\lambda_0),\bfH)\,,
$$
from which we deduce the following.
\Bt Let $(\bfw,p,\xi_\sub{\bfw}\,\bfe_1,\bfomega_\sub{\bfw})$ be a solution to the ``criticality" equation \eqref{2.7}. A sufficient condition for symmetry breaking bifurcation to occur is that
$$
\xi_0\,( \partial_1\bfw,\bfH)+2(\bfw\cdot\mathbb D(\bfH),\bfv_0)+\xi_\sub{\bfw}(\partial_1\bfv_0(\lambda_0),\bfH)\neq 0\,,
$$
where $\bfH$ is given in \eqref{A0}.
\ET{2.2}
  
\ed

\Bl Let ${\sf s}_0:=(\bfv_0,p_0,\bfxi_0,\0)$ be the solution in \theoref{1} corresponding to $\lambda=\lambda_0$. Then $\bfv_0\in\mathscr X(\Omega)$. 
\EL{3}
{\em Proof.} We only need to show
\be
\|\partial_1\bfv_0\|_{-1}<\infty\,.
\eeq{9}
We recall that ${\sf s}_0$ solves the following problem
\be\ba{cc}\medskip\left.\ba{ll}\medskip\Delta\bfv_0+\lambda_0\,\xi_0\,\partial_1\bfv_0=\lambda_0\,\bfv_0\cdot\nabla\bfv_0+\nabla p_0\\
\Div\bfv_0=0\ea\right\}\ \ \mbox{in $\Omega$}\\ \medskip
\bfv_0=\xi_0\,\bfe_1\ \ \mbox{at $\partial\Omega$}\,,\ \ \Lim{\mbox{\footnotesize $|\bfx|$}\to\infty}\bfv_0(x)=\0\,,\\ 
\Int{\partial\Omega}{}\mathbb T(\bfv_0,p_0)\cdot\bfn=\bfe_1\,,\ \   
\Int{\partial\Omega}{}\bfx\times\mathbb T(\bfv_0,p_0)\cdot\bfn=0\,.
\ea
\eeq{10}
If we dot-multiply both sides of \eqref{10}$_1$ by $\bfphi\in\calc$, integrate by parts over $\Omega$ and use \eqref{10}$_{2,3,5,6}$ we get
$$
\lambda_0(\partial_1\bfv_0,\bfphi)=(\mathbb D(\bfv_0),\mathbb D(\bfphi))-\bfphi_1\cdot\bfe_1-\lambda_0(\bfv_0\cdot\mathbb D(\bfphi),\bfv_0)\,.
$$
Thus, using  in the latter Schwarz inequality and \lemmref{1},  we infer
$$
\lambda_0|(\partial_1\bfv_0,\bfphi)|\le c\,\left(\|\mathbb D(\bfv_0)\|_2+1+\lambda_0\,\|\bfv_0\|_4^2\right)\|\mathbb D(\bfphi)\|_2\,,
$$
and since, by \theoref{1}, $\bfv_0\in L^4(\Omega)$, the proof of the lemma is completed. 
\par\hfill$\square$\par
\Bl

\EL{}